\documentclass[12pt]{article}
\usepackage{amsmath,amsfonts,amssymb,amscd}
\usepackage{graphicx}

\newtheorem{theo}{Theorem}
\newtheorem{lem}{Lemma}[section]

\newtheorem{cor}{Corollary}[section]
\newtheorem{rem}{Remark}[section]
\newtheorem{dfn}{Definition}[section]

\makeatletter \@addtoreset{equation}{section} \makeatother

\newcommand{\mC}{\mathbb{C}}
\newcommand{\mL}{\mathbb{L}}

\newcommand{\mR}{\mathbb{R}}

\newcommand{\mZ}{\mathbb{Z}}
\newcommand{\mN}{\mathbb{N}}

\newcommand{\bF}{{\bf F}}

\newcommand{\bH}{{\bf H}}

\newcommand{\bV}{{\bf V}}

\newcommand{\be}{{\bf e}}

\newcommand{\bh}{{\bf h}}

\newcommand{\bk}{{\bf k}}
\newcommand{\bl}{{\bf l}}
\newcommand{\bm}{{\bf m}}

\newcommand{\bs}{{\bf s}}

\newcommand{\bv}{{\bf v}}
\newcommand{\bw}{{\bf w}}

\newcommand{\bz}{{\bf z}}

\newcommand{\calA}{{\cal A}}

\newcommand{\calF}{{\cal F}}

\newcommand{\calH}{{\cal H}}

\newcommand{\calL}{{\cal L}}

\newcommand{\calN}{{\cal N}}

\newcommand{\calZ}{{\cal Z}}

\newcommand{\eps}{\varepsilon}
\newcommand{\ph}{\varphi}
\newcommand{\thet}{\vartheta}

\newcommand{\one}{{\bf{1}}}

\newcommand{\rank}{\operatorname{rank}}

\newcommand{\sign}{\operatorname{sign}}

\newcommand\qed{{\unskip\nobreak\hfil\penalty50
  \hskip2em\hbox{}\nobreak\hfil\mbox{\rule{1ex}{1ex} \qquad}
    \parfillskip=0pt \finalhyphendemerits=0\par\medskip}}

\begin{document}

\title
{Normalization flow in the presence of a resonance}
\author{ Dmitry Treschev \\
Steklov Mathematical Institute of Russian Academy of Sciences
}
\date{}
\maketitle

\begin{abstract}
Following \cite{Tre23}, we develop an approach to the Hamiltonian theory of normal forms based on continuous averaging. We concentrate on the case of normal forms near an elliptic singular point, but unlike \cite{Tre23} we do not assume that frequences of the linearized system are nonresonant. We study analytic properties of the normalization procedure. In particular we show that in the case of a codimension one resonance an analytic Hamiltonian function may be reduced to a normal form up to an exponentially small reminder with explicit estimates of the reminder and the analyticity domain.
\end{abstract}

\section{Introduction}

Consider a Hamiltonian system with $n$ degrees of freedom in a neighborhood of an elliptic singular point. In the linear approximation the dynamics is determined by the Hamiltonian equations
\begin{eqnarray*}
& \dot x = \partial_y H_2, \quad
  \dot y = - \partial_x H_2, \qquad
  x = (x_1,\ldots,x_n), \quad y = (y_1,\ldots,y_n), & \\
& H_2 = \sum_{j=1}^n \frac{\omega_j}{2} (y_j^2 + x_j^2). &  
\end{eqnarray*}
The real numbers $\omega_1,\ldots,\omega_n$ form the frequency vector $\omega\in\mR^n$.

According to Birkhoff \cite{Birk} it is convenient to use complex coordinates
$$
  (z,\overline z) = (z_1,\ldots,z_n,\overline z_1,\ldots,\overline z_n) \in \mC^{2n}, \qquad
  z_j = \frac{y_j+ix_j}{\sqrt 2}, \quad
  \overline z_j = \frac{y_j-ix_j}{\sqrt 2}.
$$
Below all the phase variables including $x$ and $y$ may be complex. Hence the overbar does not mean necessarily complex conjugation. In the other words, the variables $z$ and $\overline z$ are assumed to be independent.

The variables $(z,\overline z)$ are canonical i.e., for any two smooth functions
$F = F(z,\overline z), G = G(z,\overline z)$ the Poisson bracket has the form
$$
    \{F,G\}
  = i\sum_{j=1}^n
     \big(\partial_{\overline z_j} F \partial_{z_j} G - \partial_{z_j} F\partial_{\overline z_j} G
     \big).
$$
The same equation determines the bracket $\{\,,\}$ on the space of formal power series in $z$ and $\overline z$. The function $H_2$ takes the form
$$
  H_2 = \sum_{j=1}^n \omega_j z_j \overline z_j.
$$
The Hamiltonian function is assumed to be $H_2 + \widehat H$,
$$
  \widehat H= \sum_{|k|+|\overline k|\ge 3} H_{k,\overline k} z^k\overline z^{\overline k}.
$$
Here $k,\overline k\in\mZ_+^n$, $\mZ_+ :=\{0,1,\ldots\}$ are multiindices and $|k| = |k_1|+\ldots+|k_n|$ is the $l^1$-norm. Below we use the shorter notation
\begin{equation}
\label{Hdiamond}
  \bz = (z,\overline z)\in\mC^{2n}, \quad
  \bk = (k,\overline k)\in\mZ_+^{2n}, \quad
  \bz^\bk = z^k\overline z^{\overline k}, \qquad
  \widehat H = \sum_{|\bk|\ge 3} \widehat H_\bk \bz^\bk .
\end{equation}

The Hamiltonian equations are as follows:
\begin{equation}
\label{ham_eq}
  \dot z = i\partial_{\overline z} (H_2 + \widehat H), \quad
  \dot{\overline z} = - i\partial_z (H_2 + \widehat H).
\end{equation}
More generally, $\dot F = \{\widehat H,F\}$ for any function (or formal power series) $F = F(\bz)$.

According to the theory of normal forms the function $\widehat H$ may be simplified by passage to another coordinate system, \cite{Birk}.

The monomial $\widehat H_\bk\bz^\bk$ in the expansion (\ref{Hdiamond}) is said to be resonant if $\langle\omega,\overline k - k\rangle = 0$, where $\langle\,,\rangle$ is the standard inner product in $\mR^n$. The integer number $|\overline k - k|$ is said to be the order of the resonance. Any integer vector $\bk=(k,\overline k)$ which determines a resonant monomial $\bz^\bk$ satisfies
$\overline k - k \in\mL_\omega$, where
\begin{equation}
\label{mL}
  \mL_\omega = \{q\in\mZ^n : \langle\omega,q\rangle = 0 \}.
\end{equation}
If $\langle\omega,q\rangle = 0$, $q\in\mZ^n$ implies $q=0$, the frequency vector $\omega$ is called nonresonant.

By using an appropriate coordinate change one gets rid of any finite set of nonresonant monomials. To eliminate all nonresonant monomials, we have to use coordinate changes in the form of, in general, divergent power series. The Hamiltonian function $N$ obtained as a result of this (formal) normalization is called the normal form of the original Hamiltonian $H_2+\widehat H$. The normal form is unique as a formal series.

The problem of convergence/divergence of the normalizing transformation under the assumption of analyticity of $\widehat H$ is central in the theory. H.~Eliasson attracted attention of specialists to another (harder) problem: convergence/divergence of the normal form. If the normalization converges and the lattice $\calL_\omega$ is at most 1-dimensional then the system is locally completely integrable. Various versions of the inverse statement are proved in \cite{Vey,Ito,El,KKN}.

Another corollary from convergence of the normalization is Lyapunov stability of the equilibrium position in the case of a nonresonant frequency vector\footnote
{under some simple explicit additional conditions}
. Papers \cite{Fayad,Koz} contain examples of real-analytic Hamiltonians $H_2 + \widehat H$ such that the origin is Lyapunov unstable in the system (\ref{ham_eq}) although in the linear approximation the system is obviously stable.

Convergence of the normal form does not imply convergence of the normalization. But it has interesting dynamical consequences: the measure of the set covered by KAM-tori turns out to be noticeably bigger than in the case when the normal form diverges \cite{Krik}.

Convergence of the normalizing transformation is an exceptional phenomenon in any reasonable sense. At the moment this exceptionality is known in terms of Baire category \cite{Sie} and $\Gamma$-capacity \cite{PM,Krik}. Explicit examples of real-analytic systems with divergent normal form can be found in \cite{Gong,Yin,Fayad}.

The case of the ``trivial'' normal form $N = H_2$ is special. According to the Bruno-Russmann theorem \cite{Bryuno,Russ}, see also \cite{Stol}, if $N = H_2$ and $\omega$ satisfies some (rather weak) Diophantine conditions,\footnote{
usually called the Bruno conditions}
then the normalization converges.

The normalizing change of variables is traditionally constructed as a composition of an infinite sequence of coordinate changes which normalize the Hamiltonian function up to a remainder of a higher and higher degree \cite{Birk}. In further works (see for example \cite{Si-Mo}) this change of coordinates is represented as a formal series in $z$ and $\overline z$. We propose another approach to the normalization.
\smallskip

Let $\calF$ be the space of all power series in the variables $z$ and $\overline z$. Sometimes we refer to elements of $\calF$ as functions although they are only formal power series. Let $\calF_\diamond\subset\calF$ be the subspace of series which start from terms of degree at least 3.
In this paper, assuming that $H_0$ is fixed, we study various normalization flows $\phi_\xi^\delta$ on the space $\calF_\diamond$, $\delta\in [0,+\infty)$. Any shift
$$
  H_2 + \widehat H \mapsto H_2 + \phi_\xi^\delta(\widehat H), \qquad
  \widehat H\in\calF_\diamond
$$
is a transformation of the Hamiltonian function $H_2 + \widehat H$ according to a certain (depending on $\delta$ and on the initial Hamiltonian $H_2 + \widehat H$) canonical change of variables.

Any flow $\phi_\xi^\delta$ is determined by a certain ODE in $\calF_\diamond$:
\begin{equation}
\label{main}
    \partial_\delta H
  = - \{\xi H, H_2 + H\}, \qquad
    H|_{\delta=0} = \widehat H.
\end{equation}
Here $\xi$ is a linear operator on $\calF_\diamond$. In fact, (\ref{main}) is an initial value problem (IVP) for a differential equation presented in the form of Lax L-A pair. Usually such systems are considered in the theory of integrable systems. However integrability of system (\ref{main}) whatever it means seems, is of no use to us.

Let $\calN_\diamond\subset\calF_\diamond$ be the subspace which consists of series which contain only resonant monomials. We will choose operators $\xi$ such that the space $\calN_\diamond$ is an invariant manifold with respect to $\phi^\delta_\xi$ and any point of $\calN_\diamond$ is fixed. We proposed such an approach in \cite{Tre23}. In this paper we present further results. In particular here we do not assume that $\omega$ is a non-resonant vector.

In Section \ref{sec:xi} we define the operator $\xi=\xi_*$. We represent the system (\ref{main}) in the form of an infinite ODE system (\ref{aver3}) for the coefficients $H_\bk$. A more convenient equivalent form of it is the system (\ref{aver4}). Because of a special ``nilpotent'' structure of the system (\ref{aver4}), the existence and uniqueness of a solution for the corresponding IVP for any initial condition turns out to be a simple fact (Section \ref{sec:formal}). In the product topology the solution tends to the normal form as $\delta\to +\infty$.

We are particularly interested in the restriction of $\phi^\delta_{\xi_*}$ to the subspace $\calA\subset\calF$ of analytic Hamiltonian functions. In Section \ref{sec:analy} we prove that for any
$\widehat H\in\calA$ the solution $H=H(\bz,\delta)$ also lies in $\calA$ for any $\delta\ge 0$. However the polydisk of analyticity generically shrinks when $\delta$ grows. A rough lower estimate for its (poly)radius gives a quantity of order $1/(1+\delta)$ (Theorem \ref{theo:analytic}).

In Section \ref{sec:degen} we assume that the normal form $H_2 + N_\diamond$ of the Hamiltonian $H_2+\widehat H$ satisfies the equation $N_\diamond = O_r(\bz)$, $r\ge 3$. For example, generically\footnote{
if the frequency vector $\omega$ does not admit resonances of order 3}
$r=4$. The case $r>4$ is not generic but may happen for some values of parameters in multiparametric families of Hamiltonians. Assuming that $\widehat H\in\calA$, $\omega$ satisfies the Bruno Condition, and $N_\diamond=O_r(\bz)$, Theorem \ref{theo:deg_nf} says that $H(\cdot,\delta)$ stays in $\calA$ and the polydisk of analyticity has radius of order (at least) $(1+\delta)^{-1/(r-2)}$, $\delta>0$. In particular, generically this radius $\sim (1+\delta)^{-1/2}$. This improves the corresponding estimate $(1+\delta)^{-1}$ in \cite{Tre23}. 

If $\omega$ is collinear to a vector with rational components (codimention one resonance) then Theorem \ref{theo:deg_nf} implies  Corollary \ref{cor:corank1}. This corollary says that there exists a change of variables which reduces in a polydisk of radius $\sim (1+\delta)^{-1/(r-2)}$ the Hamiltonian to the form $H_2 + G^0 + G^*$, where $G^0\in\calN_\diamond$ and $|G^*|\sim e^{-c\delta}$, $c>0$.

Theorem \ref{theo:BR} from Section \ref{sec:proof_th4} is an auxiliary statement used in the proof of Theorem \ref{theo:deg_nf}. However, it probably has an independent value. It has the same assumptions: $\widehat H\in\calA$, $\omega$ satisfies the Bruno Condition, $N_\diamond = O_r(\bz)$. Then there exists a change of the variables which reduces the Hamiltonian function to the form $H_2 + G$, $G = O_r(\bz)$. Theorem \ref{theo:BR} presents upper estimates for Taylor coefficients of the function $G$ in terms of analytic properties of $\widehat H$. As a corollary we obtain the Bruno-Sigel theorem about convergence of the normalization in the case of the trivial normal form (Corollary \ref{cor:BR}).

Section \ref{sec:tech} contains technical statements concerning majorants, Bruno sequences etc.

\section{Basic construction}

For any $q = (q_1,\ldots,q_n)\in\mZ^n$ let $|q|=|q_1|+\cdots+|q_n|$ be its $l^1$-norm. For any
$k,\overline k\in\mZ_+^n$ we put
$$
  \bk=(k,\overline k)\in\mZ_+^{2n},\quad
  \bk' = \overline k - k\in\mZ^n , \quad
  \bk^* = (\overline k,k), \quad
  |\bk| = |k| + |\overline k|.
$$
Then $\bk^*-\bk=(\bk',-\bk')$. In particular $\bk'=0$ if and only if $\bk=\bk^*$.

We will use the notation
$$
  \mZ_\diamond^{2n} = \{\bk\in\mZ_+^{2n} : |\bk|\ge 3\}.
$$

\subsection{Spaces}

Let $\calF$ be the vector space of all series
\begin{equation}
\label{H_*incalF}
  H = \sum_{\bk\in\mZ_+^{2n}} H_{\bk} \bz^\bk, \qquad
  H_{\bk}\in\mC.
\end{equation}
The series (\ref{H_*incalF}) are assumed to be formal i.e., there is no restriction on the values of the coefficients $H_{\bk}$. So, $\calF$ coincides with the ring
$\mC[[z_1,\ldots,z_n,\overline z_1,\ldots,\overline z_n]]$. The Poisson bracket determines on $\calF$ the structure of a Lie algebra.
Below we use on $\calF$ the product topology i.e., a sequence $H^{(1)},H^{(2)},\ldots\in\calF$ is said to be convergent if for any $\bk\in\mZ_+^{2n}$ the sequence of coefficients
$H^{(1)}_{\bk},H^{(2)}_{\bk},\ldots$ converges.

For any $H$ satisfying (\ref{H_*incalF}) we define $p_{\bk}(H) = H_{\bk}$. So $p_{\bk} : \calF\to\mC$ is a canonical projection, corresponding to $\bk\in\mZ_+^{2n}$.
Suppose $F\in\calF$ depends on a parameter $\delta\in I$, where $I\subset\mR$ is an interval. In other words, we consider a map
$$
  f : I\to\calF, \quad I\ni\delta\mapsto F(\cdot,\delta).
$$
We say that $F$ is smooth in $\delta$ if all the maps $p_\bk\circ f$ are smooth.


Let $\calF_r\subset\calF$ be the space of ``real'' series:
$$
   \calF_r
 = \{H\in\calF : \overline H_{\bk} = H_{\bk^*} \;
                                    \mbox{ for any $\bk\in\mZ_+^n$}\}.
$$

We define $\calA\subset\calF$ as the space of analytic series:
$$
   \calA
 = \{H\in\calF : \mbox{ there exist $c,a$ such that $|H_{\bk}| \le c e^{a|\bk|}$
                          for any $\bk\in\mZ_+^{2n}$}\}.
$$

The product topology may be restricted from $\calF$ to $\calA$, but, being a scale of Banach spaces, $\calA$ may be endowed with a more natural topology. We have
$\calA = \cup_{\rho>0}\calA^\rho$, where $\calA^\rho$ is a Banach space with the norm $\|\cdot\|_\rho$,
\begin{equation}
\label{Drho}
  \|H\|_\rho = \sup_{\bz\in D_\rho} |H(\bz)|, \quad
  D_\rho = \{\bz\in\mC^{2n} : |z_j|<\rho,\; |\overline z_j|<\rho,\; j=1,\ldots,n\}.
\end{equation}
Then we have:
$\calA^{\rho'}\subset\calA^\rho$, $\|\cdot\|_\rho \le \|\cdot\|_{\rho'}$ for any
$0 < \rho < \rho'$.

\begin{lem}
\label{lem:Hkk}
If $\|H\|_\rho \le c$ then
\begin{equation}
\label{Cauchy}
  |H_{\bk}| \le c\rho^{-|\bk|}.
\end{equation}
\end{lem}

The proof is standard. It follows from the Cauchy formula: for any positive $\rho_0<\rho$ and any $H\in\calA^\rho$
$$
    H_{\bk}
  = \frac1{(2\pi i)^{2n}} \oint dz_1 \ldots \oint dz_n \oint d\overline z_1 \ldots
                  \oint \frac{H(\bz)\, d\overline z_n}{\bz^{\bk+\one}},
$$
where $\one=(1,\ldots,1)\in\mZ_+^{2n}$ and the integration is performed along the circles
$$
  |z_1|=\rho_0,\ldots,\quad
  |z_n|=\rho_0,\quad
  |\overline z_1|=\rho_0,\ldots,\quad
  |\overline z_n|=\rho_0.
$$
This implies $|H_\bk|\le c\rho_0^{-|\bk|}$. Since $\rho_0<\rho$ is arbitrary, we obtain (\ref{Cauchy}). \qed

\begin{cor}
\label{cor:stronger}
Topology determined on $\calA^\rho$ by the norm $\|\cdot\|_\rho$ is stronger than the product topology induced from $\calF$ i.e., if a sequence $\{H^{(j)}\}$ converges in $(\calA^\rho, \|\cdot\|_\rho)$ then it converges in $\calF$.
\end{cor}

We define the resonant sublattice (\ref{mL}).
The number of generators in the group $(\mL_\omega,+)$ is said to be the rank of $\mL_\omega$. In particular if $\rank\mL_\omega=0$ (equivalently, $\mL_\omega=\{0\}$) then the frequency vector is called nonresonant. If $\rank\mL_\omega=1$, we have a simple resonance if $\rank\mL_\omega>1$, the resonance is multiple, if $\rank\mL_\omega=n-1$ then we have a codimension one resonance. In the latter case all trajectories of the linearized system
$$
  \dot z = i\partial_{\overline z} H_0, \quad
  \dot {\overline z} = i\partial_z H_0
$$
are periodic.

Let $\calN\subset\calF$ be the space of ``normal forms'':
$$
   \calN
 = \{H\in\calF : H_{\bk} \ne 0 \; \mbox{ implies } \bk'\in\mL_\omega\}.
$$
It is easy to check that $\calN$ is a subalgebra in $\calF$ with respect to the operations of multiplication and the Poisson bracket.

We define the subspace (the subring) $\calF_\diamond\subset\calF$ of the series
\begin{equation}
\label{H_diamondincalF}
  H = \sum_{\bk\in\mZ_\diamond^{2n}} H_{\bk} \bz^\bk, \qquad
  H_{\bk}\in\mC
\end{equation}
i.e., for any $H\in\calF_\diamond$ we have $H=O_3(\bz)$.
By definition we also put $\calN_\diamond = \calN\cap\calF_\diamond$.

\section{Normalization flow}

Consider the change of variables in the form of a shift
\begin{equation}
\label{shift}
  \bz=(z,\overline z) \mapsto \bz_\delta=(z_\delta,\overline z_\delta)
\end{equation}
along solutions of the Hamiltonian system\footnote{
The so-called, Lie method.} with Hamiltonian $F=F(\bz,\delta)=O_3(\bz)$ and independent variable $\delta$:
\begin{equation}
\label{ham=F}
            z' = i \partial_{\overline z} F, \quad
  \overline z' = - i \partial_z F, \qquad
  (\cdot)' = d/d\delta.
\end{equation}
Suppose the function $H_2+\widehat H$ expressed in the variables $\bz_\delta$ takes the form $H_2+H$:
\begin{equation}
\label{H=H}
    H_2(\bz) + \widehat H(\bz)
  = H_2(\bz_\delta) + H(\bz_\delta,\delta).
\end{equation}

Differentiating equation (\ref{H=H}) in $\delta$, we obtain:
$$
    \partial_\delta H
  = -\{F,H_2+H\}, \qquad
    H|_{\delta=0} = \widehat H.
$$

The main idea of the continuous averaging is to take $F$ in the form $\xi H$, where $\xi$ is a certain operator\footnote
{This idea is similar to the idea of Moser's homotopy method \cite{Moser}.}
, depending on the problem we deal with. Then we obtain an initial value problem in $\calF_\diamond$:
\begin{equation}
\label{aver0}
  \partial_\delta H = - \{\xi H,H_2 + H\}, \qquad
  H|_{\delta=0} = \widehat H.
\end{equation}

\subsection{Operator $\xi_*$}
\label{sec:xi}

To explain the idea, we start with the ``simplest'' operator $\xi=\xi_*$. We put
\begin{equation}
\label{sigmaomega}
  \sigma_q = \sign (\langle\omega,q\rangle), \quad
  \omega_q = |\langle\omega,q\rangle|, \qquad
  q\in\mZ^n.
\end{equation}

For any $H$ satisfying (\ref{H_diamondincalF}) we put
\begin{equation}
\label{xi}
     \xi_* H
 = - i \sum_{\bk\in\mZ^{2n}_\diamond}
                  \sigma_{\bk'} H_{\bk} \bz^\bk
 = i (H^- - H^+), \qquad
   H^\pm
 = \sum_{\pm\sigma_{\bk'} > 0} H_{\bk} \bz^\bk.
\end{equation}
We also put
$$
    H^0
 =  \sum_{\langle\omega,\bk'\rangle = 0} H_{\bk} \bz^\bk.
$$
So that $H = H^- + H^0 + H^+$.

We obtain a more detailed form of the IVP (\ref{aver0}):
\begin{eqnarray}
\label{aver1}
&   \partial_\delta H
  = v_0(H) + v_1(H) + v_2(H), \qquad
    H|_{\delta=0} = \widehat H, & \\
\nonumber
&   v_0
  = - \{\xi_* H,H_2\}, \quad
    v_1
  = - \{\xi_* H,H^0\}, \quad
    v_2
  = - \{\xi_* H,H^- + H^+\}.
\end{eqnarray}
Informal explanation for the choice (\ref{xi}) of the operator $\xi=\xi_*$ is as follows. Removing nonlinear terms $v_1+v_2$ in (\ref{aver1}), we obtain the equation
$\partial_\delta H = v_0(H)$ or, in more detail
\begin{equation}
\label{aver_trun}
    \partial_\delta H_{\bk}
  = -\omega_{\bk'} H_{\bk}, \qquad
    H_{\bk}|_{\delta=0} = \widehat H_{\bk}
\end{equation}
which can be easily solved:
$$
  H_{\bk} = e^{-\omega_{\bk'}\delta} \widehat H_{\bk} .
$$
Thus when $\delta\to +\infty$, solution $H$ of the truncated problem (\ref{aver_trun}) tends to $H^0\in\calN_\diamond$.
\smallskip

Let $e_j=(0,\ldots,1,\ldots,0)$ be the $j$-th unit vector in $\mZ_+^n$ and let $\be_j=(e_j,e_j)\in\mZ_+^{2n}$. Equation (\ref{aver1}) written for each Taylor coefficient $H_{\bk}$ has the form
\begin{eqnarray}
\label{aver3}
     \partial_\delta H_{\bk}
 &=& -\omega_{\bk'} H_{\bk} + v_{1,\bk}(H) + v_{2,\bk}(H),\qquad
     H_\bk|_{\delta=0} = \widehat H_\bk, \\
\label{v0}
     v_{1,\bk}
 &=& - \sum_{j=1}^n \sum_{\sigma_{\bm'}=0,\, \bl+\bm-\bk=\be_j} \sigma_{\bk'}
                   (\overline l_j m_j - l_j\overline m_j) H_{\bl} H_{\bm}, \\
\nonumber
     v_{2,\bk}
 &=& 2 \sum_{j=1}^n \sum_{\sigma_{\bl'} < 0 < \sigma_{\bm'},\,
                                 \bl+\bm-\bk = \be_j}
         (\overline l_j m_j - l_j\overline m_j) H_{\bl} H_{\bm}.
\end{eqnarray}

By using the change of variables
$$
  H_{\bk} = \calH_{\bk} e^{-\omega_{\bk'} \delta},
$$
we reduce the equations (\ref{aver3}) to the form
\begin{eqnarray}
\label{aver4}
\!\!\!\!\!\!\!
     \partial_\delta\calH_{\bk}
 &=&  v_{1,\bk}(\calH) + \bv_{2,\bk}(\calH), \qquad
     \calH_{\bk}|_{\delta=0} = \widehat H_{\bk}, \\
\label{bv2}\!\!\!\!\!\!\!
     \bv_{2,\bk}
 &=& 2 \sum_{j=1}^n \sum_{\sigma_{\bl'} < 0 < \sigma_{\bm'},\,
                                 \bl+\bm-\bk = \be_j}
     \!  (\overline l_j m_j - l_j\overline m_j)\calH_{\bl}\calH_{\bm}
          e^{- \omega_{\bl',\bm'} \delta}, \\
\nonumber\!\!\!\!\!\!\!
     \omega_{\bl',\bm'}
 &=& \omega_{\bl'}+\omega_{\bm'}-\omega_{\bl'+\bm'}
\end{eqnarray}

\begin{rem}
\label{rem:1-5}

1. The functions $v_{1,\bk}$ and $v_{2,\bk}$ are quadratic polynomials in the variables $H_{\bk}$. The functions $\bv_{2,\bk}$ are quadratic polynomials in $\calH_{\bk}$ with coefficients, depending on $\delta$.

2. The polynomial $v_{1,\bk}$ vanishes if $\bk'\in\mL_\omega$.

3. The polynomials $v_{1,\bk}, v_{2,\bk}$ and $\bv_{2,\bk}$ do not depend on the variables $H_\bs,\calH_{\bs}$, for any $\bs\in\mZ_\diamond^{2n}$ such that $|\bs| > |\bk|-2$.

4. The polynomials $v_{2,\bk}$ and $\bv_{2,\bk}$ do not depend on the variables $H_\bs,\calH_{\bs}$ for any $\bs\in\mZ_\diamond^{2n}$ such that
$0\le\sigma_{\bk'} \langle\omega,\bs'\rangle\le\omega_{\bk'}$.

5. For any $\sigma_{\bl'} < 0 < \sigma_{\bm'}$ the quantity $\omega_{\bl',\bm'}$ is positive:
$$
    \omega_{\bl',\bm'}
  = \left\{\begin{array}{ccc}
           2\omega_{\bl'} &\mbox{ if }& \sigma_{\bk'} > 0, \\
           2\omega_{\bm'} &\mbox{ if }& \sigma_{\bk'} < 0.
           \end{array}
    \right.
$$
\end{rem}

\section{Resonant normalization exists}

\subsection{Formal aspect}
\label{sec:formal}

Let $I\subset\mR$ be an interval containing the point $0$. We say that the path $\gamma:I\to\calF_\diamond$ is a solution of the system (\ref{aver0}) on the interval $I$ if $\gamma(0)=\widehat H$ and the coefficients
$H_{\bk}(\delta) = p_{\bk}\gamma(\delta)$ satisfy (\ref{aver3}).

If the solution $\gamma : I\to\calF_\diamond$ exists and is unique then for any $\delta\in I$ we put $\gamma(\delta)=\phi^\delta_{\xi_*}(\widehat H)$.

\begin{dfn}
\label{dfn:nil}
The ordinary differential equation
$$
  F' = \Phi(F), \qquad
  F = \sum_{|\bk|\ge 3} F_\bk \bz^\bk
$$
on $\calF$ has a nilpotent form if for any $\bk\in\mZ_\diamond^{2n}$ we have: $F'_\bk=\Phi_\bk$, where $\Phi_\bk = p_\bk\circ\Phi$ is a function of the variables $F_\bm$, $|\bm|<|\bk|$.
\end{dfn}

According to item 3 of Remark \ref{rem:1-5} the system (\ref{aver4}) has a nilpotent form. In particular, $\partial_\delta\calH_{\bk} = 0$ for $|\bk|=3$.

\begin{theo}
\label{theo:phi}
For any $\widehat H\in\calF_\diamond$ and any $\delta\in\mR$ the element $H(\cdot,\delta)=\phi^\delta_{\xi_*}(\widehat H)\in\calF$ is well-defined. For any $\bk\in\mZ_\diamond^{2n}$ the function $\calH_\bk(\delta) = H_\bk(\delta) e^{\omega_{\bk'}\delta}$ equals
\begin{equation}
\label{Pk}
  \calH_\bk(\delta) = \widehat H_\bk + P_\bk(\widehat H,\delta),
\end{equation}
where $P_\bk$ is a polynomial in $\widehat H_\bm$, $|\bm|<|\bk|$ with coefficients in the form of (finite) linear combinations of terms $\delta^s e^{-\nu\delta}$, $s\in\mZ_+$. Here $\nu\ge 0$ and moreover, if the term $\delta^s e^{-\nu\delta}$ with $\nu=0$ appears in a coefficient of $P_\bk$, $\bk'\in\mL_\omega$ then in this term $s=0$.
\end{theo}

{\it Proof of Theorem \ref{theo:phi}}. For $|\bk|=3$ we have: $\calH_\bk(\delta)=\widehat H_\bk$.
 Suppose equations (\ref{Pk}) hold for all vectors $\bk\in\mZ_\diamond^{2n}$ with $|\bk|<K$. Take any $\bk\in\mZ^{2n}_\diamond$ with $|\bk|=K$. By (\ref{aver4})
$$
  \calH_\bk=\widehat H_\bk + I_1 + I_2, \qquad
  I_1 = \int_0^\delta v_{1,\bk}(\calH(\lambda))\, d\lambda, \quad
  I_2 = \int_0^\delta
          \bv_{2,\bk}(\calH(\lambda),\lambda)\, d\lambda.
$$
By using (\ref{v0}), (\ref{bv2}), and the induction assumption, we obtain (\ref{Pk}).

If $\bk'\in\mL_\omega$ then $v_{1,\bk}=0$ while each term in $\bv_{2,\bk}$ has the form
$\lambda^s e^{-\lambda\delta}$ with $\lambda>0$. This implies the last statement of Theorem \ref{theo:phi}.
\qed

\begin{cor}
\label{cor:phi}
The limit $\lim_{\delta\to +\infty} \phi^\delta(\widehat H)$ exists in the product topology on $\calF_\diamond$ and lies in $\calN_\diamond$.
\end{cor}

Indeed, by Theorem \ref{theo:phi} for any $\bk\in\mZ_\diamond^{2n}$ there exists the limit
$$
  \lim_{\delta\to +\infty} H_{\bk}(\delta) = H_{\bk}(+\infty).
$$
The convergence is exponential and $H_{\bk}(+\infty)$ vanishes if $\bk'\not\in\mL_\omega$. This is equivalent to the existence of
$\lim_{\delta\to +\infty}\phi^\delta_{\xi_*}(\widehat H)=\phi^{+\infty}_{\xi_*}(\widehat H)$
together with the statement $\phi^{+\infty}_{\xi_*}(\widehat H)\in\calN_\diamond$.

\begin{cor}
\label{cor:real}
Suppose $\widehat H\in\calF_r\cap\calF_\diamond$. Then
$\phi^\delta_{\xi_*}(\widehat H)\in\calF_r\cap\calF_\diamond$ for any $\delta\ge 0$.
\end{cor}

Indeed, the required identities $\overline H_{\bk}(\delta) = H_{\bk^*}(\delta)$, $\delta\ge 0$ can be proved by induction in $|\bk|$ by using the equation (\ref{aver4}).

\subsection{Analytic aspect}
\label{sec:analy}

\begin{theo}
\label{theo:analytic}
Suppose $\widehat H\in\calA^\rho\cap\calF_\diamond$, $\|\widehat H\|_\rho=h\rho^3$. Then for any $\delta\ge 0$
$$
  \calH = \sum\calH_\bk(\delta)\bz^\bk\in\calA^{g(\delta)}\cap\calF_\diamond,
$$
where for some constants $A,B>0$, depending on $\rho,h$, and $n$
\begin{equation}
\label{g>A}
  g(\delta) \ge \frac A{1 + B\delta}.
\end{equation}
More precisely,
$$
      |\calH|
  \le \frac{h\rho^3}{4 (1 + 32 nh\rho\delta)^3}
  \quad\mbox{in the domain}\quad
  \Big\{\bz\in\mC^{2n} : |\bz| \le \frac{\rho}{4(1 + 32 nh\rho\delta)}
  \Big\} .
$$
\end{theo}

{\it Proof}. We use the majorant method. We remind definitions and basic facts concerning majorants in Section \ref{sec:maj}. We have:
$$
      \widehat H
  \ll f(\zeta)
   =  \sum_{\bk\in\mZ_\diamond^{2n}} \widehat\bH_{\bk} \bz^\bk, \qquad
  \zeta = \sum_{j=1}^n (z_j + \overline z_j) .
$$
The function $f(\zeta)=O_3(\zeta)$, a majorant for the initial condition, will be chosen later.

Consider together with (\ref{aver4}) the following majorant system
\begin{eqnarray}
\label{aver_maj}
\!\!\!\!\!\!
     \partial_\delta\bH_{\bk}
 &=&  \bV_{\bk}, \qquad
     \bH_{\bk}|_{\delta=0} = \widehat\bH_{\bk}, \\
\label{bV}
\!\!\!\!\!\!
     \bV_{\bk}
 &=& 2 \sum_{j=1}^n \sum_{\bl+\bm-\bk = \be_j}
         (\overline l_j m_j + l_j\overline m_j)\bH_{\bl}\bH_{\bm}.
\end{eqnarray}
To obtain equation (\ref{bV}), we have replaced in $v_{1,\bk}$ and $\bv_{2,\bk}$ minuses by pluses, dropped the exponential multipliers, and added some new positive (for positive $\bH_{\bl}$ and $\bH_{\bm}$) terms. The main property of the system (\ref{aver_maj})--(\ref{bV}) is such that if for any $\bk\in\mZ_\diamond^{2n}$ we have
$|\widehat H_{\bk}|\le\widehat\bH_{\bk}$ then for any $\delta\ge 0$ the inequality
$|v_{1,\bk} + \bv_{2,\bk}| \le \bV_\bk$ holds i.e., conditions (a) and (b) from Definition \ref{dfn:maj} hold.

The system (\ref{aver_maj})--(\ref{bV}) has a nilpotent form. Therefore by Theorem \ref{theo:maj_nil} we may use Majorant principle: if $\bH$ is a solution of (\ref{aver_maj})--(\ref{bV}) then (\ref{aver4}) has a solution $\calH$ and $\calH\ll\bH$ for any $\delta\ge 0$.

The system (\ref{aver_maj})--(\ref{bV}) may be written in a shorter form:
\begin{equation}
\label{maj_shorter}
  \partial_\delta\bH = 4\sum_{j=1}^n \partial_{z_j}\bH\, \partial_{\overline z_j}\bH, \qquad
  \bH|_{\delta=0} = f(\zeta).
\end{equation}
Since the initial condition $\bH|_{\delta=0}$ depends on $z$ and $\overline z$ only through the variable $\zeta$, we may look for a solution of (\ref{maj_shorter}) in the form
$\bH(z,\overline z,\delta)=F(\zeta,\delta)$. The function $F$ satisfies the equation
$$
  \partial_\delta F = 4n (\partial_\zeta F)^2, \qquad
  F|_{\delta=0} = f(\zeta).
$$

The function $G=\partial_\zeta F$ satisfies the inviscid Burgers' equation
\begin{equation}
\label{Burgers}
  \partial_\delta G = 8n G \partial_\zeta G, \qquad
  G|_{\delta=0} = \partial_\zeta f(\zeta) = O_2(\zeta).
\end{equation}

By using the method of characteristics we obtain that the function $G=G(x,t)$ which solves (\ref{Burgers}), satisfies the equation
\begin{equation}
\label{Bur_sol}
  G = f'(\zeta + 8n\delta G).
\end{equation}

By Lemma \ref{lem:maj1} we can take $f(\zeta) = h\rho\zeta^3 / (\rho-\zeta)$. Then by Lemma \ref{lem:f'}
$$
  f'(\zeta) \ll a\zeta^2 / (b-\zeta), \qquad
  a = 2h\rho, \quad  b= \rho/2.
$$
Putting $\tau=8n\delta$, we obtain from (\ref{Bur_sol})
$$
  G = \frac{a(\zeta+\tau G)^2}{b - \zeta - \tau G}.
$$
This is a quadratic equation w.r.t. $G$. The solution is
$$
  G = \frac{2a\zeta^2}{b-\zeta-2a\tau\zeta + \sqrt{(b-\zeta-2a\tau\zeta)^2 - 4a\tau\zeta^2(1+a\tau)}}.
$$
It is analytic for
$$
  |\zeta| < \frac{b}{1+2a\tau + 2\sqrt{a\tau(1+a\tau)}} .
$$

If $\displaystyle |\zeta| \le \frac{b}{2(1+2a\tau)}$ then
$\displaystyle |G| \le \frac{ab}{(1 + 2a\tau)^2}$. This implies that
$$
      |\calH|
  \le |F|
  \le \frac{ab^2}{2(1 + 2a\tau)^3}\quad
  \mbox{in the domain}\quad
  \Big\{|\bz| \le \frac{b}{2(1+2a\tau)} \Big\}.
$$
\qed

\begin{cor}
\label{cor:ra}
Suppose $\widehat H\in\calA\cap\calF_r$. Then
$\phi^\delta_{\xi_*}(\widehat H)\in\calA\cap\calF_r$ for any $\delta\ge 0$.
\end{cor}

Note that for a typical $\widehat H\in\calA\cap\calF_\diamond$ one should expect that
$\phi^{+\infty}_{\xi_*}(\widehat H)\in\calN_\diamond$ does not belong to $\calA$, \cite{Krik}.

\section{Degenerate normal form}
\label{sec:degen}

\subsection{Normalization theorem}

Suppose the normal form $H_2 + N_\diamond$ is such that
$N_\diamond := \lim_{\delta\to +\infty} \phi_{\xi_*}^\delta (\widehat H) \in \calN$ vanishes in several first orders in $\bz$. In this section we obtain more precise estimates (in comparison with those from Theorem \ref{theo:analytic}) for the analyticity domain of the Hamiltonian function $\calH(\cdot,\delta)$ when $\delta$ increases.

For any $j\in\mN$ we put
\begin{equation}
\label{Omega}
    \Omega_s
  = \max \Big\{ \frac1{|\langle\omega,q\rangle|} :
                  q\in\mZ^n\setminus\mL_\omega, \, |q|\le s \Big\}.
\end{equation}

\begin{dfn}
We say that $\{a_j\}_{j\in\mZ_+}$, $a_j\ge 1$ is a Bruno sequence if $\{a_j\}$ is monotone non-decreasing and
$$
  \sum_{j=0}^\infty 2^{-j} \ln a_j < \infty.
$$
\end{dfn}

Obviously if $\{a_j\}$ and $\{b_j\}$ are two Bruno sequences then the product $\{a_j b_j\}$ is also a Bruno sequence.

\begin{dfn}
We say that the vector $\omega\in\mR^n$ satisfies the Bruno condition if the sequence $\{a_j\}$, $a_j=\max\{1,\Omega_{2^j+2}\}$ is Bruno.
\end{dfn}

\begin{theo}
\label{theo:deg_nf}
Suppose $\widehat H\in\calA^\rho\cap\calF_\diamond$, $\omega$ satisfies the Bruno condition, and
$N_\diamond = O_r(\bz)$. Then for any $\delta\ge 0$ the solution
$\calH(\bz,\delta) =  \sum \calH_\bk(\delta) \bz^\bk$  of (\ref{aver4}) belongs to the space
$\calA^{g(\delta)}$, where for some constant $C$, depending on $\rho,\|\widehat H\|_\rho,\{\Omega_j\}$, and $n$
\begin{equation}
\label{g>C/delta}
  g(\delta) \ge \frac{C}{1 + \delta^{1/(r-2)}}.
\end{equation}
More precisely, for some $C_H = C_H(\rho,\|\widehat H\|_\rho,\{\Omega_j\},n)$ and
$C_z = C_z(\rho,\|\widehat H\|_\rho,\{\Omega_j\},n)$
\begin{equation}
\label{|H|<}
  |\calH| \le \frac{C_H}{1 + \delta^{r/(r-2)}} \quad
  \mbox{in the domain}\quad
  \Big\{ \bz\in\mC^{2n} : |\bz| \le \frac{C_z}{1+\delta^{1/(r-2)}} \Big\}.
\end{equation}
\end{theo}

In particular, if $r=3$, Theorem \ref{theo:deg_nf} essentially coincides with Theorem \ref{theo:analytic}. If $r\ge 4$, (\ref{g>C/delta}) improves the estimate (\ref{g>A}).

Proof of Theorem \ref{theo:deg_nf} is presented in Section \ref{sec:proof_th4}.

\subsection{Resonance of codimension one}

Suppose $\rank\mL_\omega = n-1$. Then there exists a unique collection of integers
$q_1,\ldots,q_n,p$, and $\lambda>0$ such that
$$
  p>0, \quad
  GCD(|q_1|,\ldots,|q_n|,p) = 1, \quad
  \mbox{and}\quad
  \omega = \lambda q / p.
$$
Hence in (\ref{bv2}) we have $\omega_{\bl',\bm'} \ge \lambda / p$.

If the normal form of $H_2+\widehat H$ is of order $r$ in $\bz$, we may use Theorem \ref{theo:deg_nf} to obtain the following corollary.

\begin{cor}
\label{cor:corank1}
Suppose $\rank\mL_\omega = n-1$, $\widehat H\in\calA^\rho\cap\calF_\diamond$,
and the normal form of $H_2+\widehat H$ is of order $O_r(\bz)$.

Then there exists a canonical change of variables $\bw\mapsto\bz$ which transforms the Hamiltonian to $H_2(\bw)+G(\bw)$, $G=O_r(\bw)$, where
$$
  G = G^0 + G^*,  \qquad
  G^0\in\calN_\diamond
$$
and in the domain
\begin{equation}
\label{domain}
  \displaystyle \Big\{\bw\in\mC^{2n} :
                           |\bw| \le \frac{C_z}{1 + \delta^{1/(r-2)}}
                \Big\}
\end{equation}
the following estimate holds:
\begin{equation}
\label{H-H0<}
  |G^0| \le \frac{C_H}{1 + \delta^{r/(r-2)}}  \quad
  \mbox{and}\quad
  |G^*| \le e^{-\lambda\delta/p} \frac{C_H}{1 + \delta^{r/(r-2)}}.
\end{equation}
\end{cor}

Estimates (\ref{H-H0<}) mean that the nonresonant part $G^*$ may be reduced to an exponentially small quantity of order $e^{-\lambda\delta/p}$. The price is a little size ($\sim \delta^{-1/(r-2)}$) of the domain in which the Hamiltonian function $G$ is defined.

To obtain estimates (\ref{H-H0<}) in the domain (\ref{domain}), we split the function $\calH(\bz,\delta)$, obtained in Theorem \ref{theo:deg_nf} into the resonant, $\calH^0$, and non-resonant, $\calH^*$, parts
$$
     \calH^0
  =  \sum_{|\bk|\ge 3,\, \bk'\in\mL_\omega} \calH_\bk \bz^\bk, \quad 
     \calH^*
  =  \sum_{|\bk|\ge 3,\, \bk'\in\mZ^n\setminus\mL_\omega} \calH_\bk \bz^\bk .
$$
The corresponding function $H(\bz,\delta)$ equals $G^0 + G^*$,
$$
     G^0
  =  \calH^0, \quad
     G^*
  =  \sum_{|\bk|\ge 3,\, \bk'\in\mZ^n\setminus\mL_\omega} \calH_\bk e^{-\omega_{\bk'}\delta} \bz^\bk .
$$
Estimates (\ref{H-H0<}) follow from (\ref{|H|<}) and the inequality $\omega_{\bk'}\ge\lambda/p$ for any $\bk'\in\mZ^n\setminus\mL_\omega$.

\section{Proof of Theorem \ref{theo:deg_nf}}
\label{sec:proof_th4}

\subsection{Inductive step}
\label{sec:induc}

\begin{dfn}
The sequence $\{b_j\}_{j\in\mN}$ is said to be sublinear if $\lim_{j\to+\infty} b_j / j = 0$.
The sequence $\{b_j\}_{j\in\mN}$ is said to be convex if $b_{j-1} - 2b_j + b_{j+1} \ge 0$ for all
$j\ge 2$.
\end{dfn}

\begin{lem}
\label{lem:low}
Suppose $\widehat H\in\calA\cap\calF_\diamond$ is such that

$(1)$ $\displaystyle\widehat H = \sum_{|\bk|\ge s} \widehat H_\bk\bz^\bk = O_s(\bz)$, \quad
        $|\widehat H_\bk| \le c e^{b_{|\bk|} + \alpha |\bk|}$, \quad
        $\alpha\ge 0$, \quad
        $s\ge 3$,

$(2)$ the sequence $\{b_j\}$ is nonincreasing, sublinear and convex\footnote{
We will use sequences $\{b_j\}$ with $b_j < 0$}.

Then an analytic canonical transformation $\bw\mapsto\bz=\thet(\bw)$ of coordinates reduces the Hamiltonian $H_2 + \widehat H$ to the form $H_2\circ\thet+\widehat H\circ\thet = H_2 + G^0 + G$, where
\begin{eqnarray}
\label{G0}
&\displaystyle
    G^0
  = \sum_{s\le |\bk| \le 2s-3,\,\langle\omega,\bk'\rangle=0}
       \widehat H_\bk \bw^\bk \in \calN_\diamond , &   \\
\label{GGk}
&\displaystyle
   G = \sum_{|\bk|\ge 2s-2} G_\bk\bw^\bk, \qquad
               |G_\bk| \le c e^{b_{|\bk|} + (\alpha+\eps) |\bk|}, & \\
\label{Lambda}
&\displaystyle
      \eps
  =   c\Lambda\Omega_{2s-2}, \qquad
      \Lambda
   =  \frac{1}{2} n e^{2\alpha} (2s)^{2n+1} e^{2b_s - b_{2s-2}}, &
\end{eqnarray}
Moreover, for any $0 < \rho\le e^{-\alpha}$ and $\rho'$ satisfying
\begin{equation}
\label{rhorho'}
  0 < \rho' \le \rho - \frac{\eps}{n s e^\alpha} (\rho e^{\alpha})^{s-1}
\end{equation}
we have: $\thet: D_{\rho'}\to D_\rho$.
\end{lem}

{\it Proof}. We construct the transformation $\thet$ by using continuous averaging with the following operator $\xi_s$. Suppose $H_2(\bz_\delta) + H(\bz_\delta,\delta)$,
$$
  H(\bz,\delta) = \sum_{|\bk|\ge s} H_\bk(\delta) \bz^\bk
$$
is the Hamiltonian function obtained as a result of the coordinate transformation $\bz_\delta\mapsto\bz$. We define
$$
  \xi_s H(\bz,\delta) = -i \sum_{\bk\in\calZ_s} \sigma_{\bk'} H_\bk(\delta) \bz^\bk, \qquad
  \calZ_s = \{\bk\in\mZ_\diamond^{2n} : s\le |\bk|\le 2s-3\},
$$
where $\sigma_{\bk'}$ is determined by (\ref{sigmaomega}).
We have: $H = H^- + H^+ + H^0 + H^*$:
\begin{eqnarray*}
& H^\pm = \sum_{\bk\in\calZ_s^\pm} H_\bk(\delta) \bz^\bk, \quad
  H^0 = \sum_{\bk\in\calZ_s^0} H_\bk(\delta) \bz^\bk, \quad
  H^* = \sum_{|\bk|\ge 2s-2} H_\bk(\delta) \bz^\bk, & \\
& \calZ_s^\pm = \{\bk\in\calZ_s : \pm\sigma_{\bk'} > 0\}, \quad
   \calZ_s^0 = \{\bk\in\calZ_s : \sigma_{\bk'} = 0\}.
\end{eqnarray*}
Then by definition $\xi_s H = i(H^- - H^+)$. Equation (\ref{aver0}) takes the form
\begin{eqnarray}
\label{pm0}
& \partial_\delta H^\pm = \pm i \{H^\pm, H_2\}, \quad
  \partial_\delta H^0 = 0, & \\
\label{*}
& \partial_\delta H^* = - i \{H^- - H^+, H^- + H^+ + H^0 + H^*\}. &
\end{eqnarray}
Here we used the fact that for any two monomials of degree at least $s$ their Poisson bracket has degree at least $2s-2$.

We present system (\ref{pm0})--(\ref{*}) in a more detailed form:
\begin{eqnarray}
\label{pm0+}
      \partial_\delta H_\bk
  &=& - \omega_{\bk'} H_\bk, \qquad\qquad\qquad\qquad\qquad\qquad\quad\;
      \mbox{for} \quad  \bk\in\calZ_s,   \\
\label{*+}
      \partial_\delta H_\bk
  &=& - i \{H^- - H^+, H^- + H^+ + H^0 + H^*\}_\bk \quad
      \mbox{for}\quad  |\bk| \ge 2s-2, \\
\label{H|delta=0}
      H_\bk|_{\delta=0}
  &=& \widehat H_\bk .
\end{eqnarray}
Here $\omega_{\bk'}$ is determined by (\ref{sigmaomega}), and $\{\,,\}_\bk$ denotes the Taylor coefficient with the number $\bk$ in the expansion of the bracket. We solve the IVP (\ref{pm0+})--(\ref{H|delta=0}) on the interval $\delta\in [0,+\infty)$. The required change of variables is obtained when $\delta=+\infty$.

Equations (\ref{pm0+}) can be easily solved: $H_\bk = e^{-\omega_{\bk'}\delta} \widehat H_\bk$ for $\bk\in\calZ_s$. This implies (\ref{G0}).

Now we associate with system (\ref{pm0+})--(\ref{*+}) a majorant system for the functions
$$
  \bH^\pm = \sum_{\bk\in\calZ_s^\pm} \bH_\bk \bz^\bk \gg H^\pm, \quad
  \bH^0 = \sum_{\bk\in\calZ_s^0} \bH_\bk \bz^\bk \gg H^0 \equiv \widehat H^0, \quad
  \bH^* = \sum_{|\bk|\ge 2s-2} \bH_\bk \bz^\bk \gg H^*.
$$
By condition (1) of Lemma \ref{lem:low} we may choose
\begin{eqnarray*}
     \widehat\bH_\bk
 &=& \bH_\bk |_{\delta=0}
  =  c e^{b_{|\bk|} + \alpha |\bk|}, \qquad
     |\bk| \ge s, \\
     \bH_\bk(\delta)
 &=& e^{-\delta/\Omega_{2s-2}}
     c e^{b_{|\bk|} + \alpha |\bk|}, \qquad\;
     \bk\in\calZ_s^- \cup \calZ_s^+ , \quad \delta\ge 0.
\end{eqnarray*}

We take as a majorant system for (\ref{*+}) the following one:
\begin{eqnarray}
\label{SigSig}
  &&   \partial_\delta\bH_\bk
 \ge  e^{-\delta / \Omega_{2s-2}} \Sigma_1(\bk), \qquad   |\bk| \ge 2s-2, \\
\nonumber
       \Sigma_1(\bk)
 &=&   \sum_{j=1}^n \sum_{\bl\in\calZ_s^-\cup\calZ_s^+,\, |\bm|\ge s,\, \bl+\bm=\bk+\be_j}
                       (\overline l_j m_j + l_j\overline m_j) \widehat\bH_\bl \bH_\bm.
\end{eqnarray}

For any $|\bk|\ge 2s-2$ we put
$\;\bH_\bk = c\, e^{b_{|\bk|} + \alpha |\bk|} \bh_\bk, \quad
 \bh_\bk(\delta) \ge \bh_\bk(0) = 1$.
\smallskip

We use the notation
$$
  \chi_K(\delta) = \max_{|\bk|\le K} \bh_\bk(\delta), \qquad
  \delta\ge 0.
$$

\begin{lem}
\label{lem:Lambda}
Suppose the sequence $\{b_j\}_{j\in\mN}$ is nonincreasing and convex. Then for any $\bk\in\mZ_+^{2n}$, $|\bk|\ge 2s-2$ the estimate
\begin{equation}
\label{Sigma1}
  \Sigma_1(\bk) \le c^2 |\bk| e^{\alpha |\bk|} \Lambda \chi_{|\bk|}
\end{equation}
holds, where $\Lambda$ is determined by (\ref{Lambda}).
\end{lem}

We prove Lemma \ref{lem:Lambda} in the end of this section.

By Lemma \ref{lem:Lambda} the differential inequalities (\ref{SigSig}) follow from (and may be replaced by) the system
$$
    \partial_\delta \bh_\bk
  =  c |\bk| e^{-\delta / \Omega_{2s-2}} \Lambda \chi_{|\bk|}, \qquad
  \bh_\bk(0) = 1.
$$
Since the multiplier $c e^{-\delta / \Omega_{2s-2}} \Lambda$ does not depend on $\bk$, we conclude\footnote{for example, by induction in $|\bk|$} that $\chi_{|\bk|} = \bh_\bk$ and we may take
$$
    \bh_\bk(\delta)
  = \exp\bigg( c |\bk| \Lambda
               \int_0^\delta e^{-\widetilde\delta / \Omega_{2s-2}}\, d\widetilde\delta \bigg).
$$
In particular, $\bh_\bk = \bh_\bk(\delta)$ is an increasing function of $\delta\in [0,+\infty)$. Taking $\delta = +\infty$, we obtain:
$$
    \bh_\bk(+\infty)
  = \exp\big( c \Lambda \Omega_{2s-2} |\bk| \big).
$$
This implies that we may use in (\ref{GGk}) $\eps$, determined by (\ref{Lambda}).
\smallskip

Let $\thet^{-1}_\delta$ be the $\delta$-shift along solutions of the Hamiltonian system
$$
  \dot z_j = i\partial_{\overline z} \xi_s H, \quad
  \dot{\overline z_j} = -i\partial_z \xi_s H, \qquad
  \xi_s H = - i \sum_{\bk\in\calZ_s} \sigma_{\bk'} \widehat H_\bk e^{-\omega_{\bk'}\delta} \bz^\bk.
$$
We estimate $\rho'$ such that for any $\bz\in D_{\rho'}$ and $\delta\in [0,+\infty)$ the point $\bz(\delta) = \thet^{-1}_\delta(\bz)$ does not leave $D_\rho$.
While $\bz(\delta)=(z_1(\delta),\ldots,\overline z_n) \in D_\rho$, $\rho e^\alpha \le 1$, we have the estimate
\begin{equation}
\label{|z|}
      |\dot z_j(\delta)|
  \le \sum_{\bk\in\calZ_s} |\bk| c e^{b_{|\bk|}+\alpha |\bk| - \delta / \Omega_{2s-2}} \rho^{|\bk|-1}
  \le (2s)^{2n} c e^{b_s - \delta/\Omega_{2s-2}} (\rho e^\alpha)^{s-1} e^\alpha,
\end{equation}
where $j$ is any number from the set $\{1,\ldots,n\}$.
Here we have used that the sequence $b_s$ does not increase, $s\le |\bk|\le 2s$, and the sum contains less than $(2s)^{2n-1}$ terms. Estimate (\ref{|z|}) implies
$$
      |z_j(\delta) - z_j(0)|
  \le (2s)^{2n} c e^{b_s} \Omega_{2s-2} (\rho e^\alpha)^{s-1} e^\alpha,  \qquad  j = 1,\ldots,n.
$$
Since the sequence $\{b_j\}$ does not increase, we have: $b_s - b_{2s-2}\ge 0$. Hence
$$
      |z_j(\delta) - z_j(0)|
  \le (2s)^{2n} c e^{2b_s - b_{2s-2}} \Omega_{2s-2} (\rho e^\alpha)^{s-1} e^\alpha
   =  \frac{1}{n s e^\alpha} (\rho e^\alpha)^{s-1} \eps .
$$
The differences $\overline z_j(\delta) - \overline z_j(0)$ can be estimated in the same way. This implies (\ref{rhorho'}). \qed

{\it Proof of Lemma \ref{lem:Lambda}}. For any $\bl,\bm,\bk\in\mZ_+^{2n}$ the equation $\bl+\bm=\bk+\be_j$ implies $|\bl|+|\bm|=|\bk|+2$. This implies the estimate
\begin{eqnarray}
   && \Sigma_1(\bk)
  \le e^{\alpha |\bk| + 2\alpha} \chi_{|\bk|} \Sigma_2(\bk), \\
\nonumber
       \Sigma_2(\bk)
 &=&   \sum_{j=1}^n \sum_{\bl\in\calZ_s^-\cup\calZ_s^+,\, |\bm|\ge s,\, \bl+\bm=\bk+\be_j}
                       c^2\, |\bl|\cdot |\bm|
                        e^{b_{|\bl|} + b_{|\bm|} - b_{|\bk|} }.
\end{eqnarray}

Since in $\Sigma_2$ $|\bm|\le |\bk|$, we have: $\Sigma_2(\bk)\le c^2 |\bk| \Sigma_3(\bk)$, where
$$
      \Sigma_3(\bk)
  =   \sum_{j=1}^n S_j(\bk), \qquad
      S_j(\bk)
  =   \sum_{\bl\in\calZ_s^-\cup\calZ_s^+,\, |\bm|\ge s,\, \bl+\bm=\bk+\be_j}
                       |\bl|\cdot e^{b_{|\bl|} + b_{|\bm|} - b_{|\bk|}}.
$$

The number of terms in the sum $S_j(\bk)$ with fixed $|\bl|\le 2s-3$ does not exceed
$$
    \# \{\bl'\in\mZ_+^{2n} : |\bl'| = |\bl|\}
  < (2s)^{2n-1}.
$$
Therefore
$$
      S_j(\bk)
  \le \sum_{p+q=|\bk|+2,\, s\le p \le 2s - 3,\, q\ge s}
                     (2s)^{2n} e^{b_p + b_q - b_{|\bk|}}.
$$
The sequence $\{b_j\}$ is convex. Therefore by Lemma \ref{lem:conv_seq}, if $|\bk| \ge 2s-2$, $p+q=|\bk|+2$, and $p,q\ge s$, we have:
$$
      b_p + b_q \le b_{|\bk|+2-s} + b_s \quad\mbox{and}\quad
      b_{|\bk|+2-s} - b_{|\bk|} \le b_s - b_{2s-2}.
$$
This implies $b_p + b_q - b_{|\bk|} \le 2b_s - b_{2s-2}$. Hence
$$
      S_j(\bk)
  \le \frac12 (2s)^{2n+1} e^{2b_s-b_{2s}-2}
   =  \frac{\Lambda e^{-2\alpha}}{n},      \quad
      \Sigma_3(\bk)
  \le e^{-2\alpha}\Lambda
$$
and we obtain (\ref{Sigma1}).  \qed

\subsection{Low order normalization}

In this section, assuming that the normal form of $H_2 + \widehat H$ is $H_2+N_\diamond$, $N_\diamond=O_r(\bz)$, we apply Lemma \ref{lem:low} inductively to reduce the Hamiltonian to the form $H_2 + G$, $G=O_r$. Assuming $\widehat H\in\calA$, we estimate $G$ in terms of $\widehat H$. As a by-product we obtain another proof of the Bruno-Russmann theorem on convergence of the normalization in the case of zero normal form.

Motivated by conditions of Theorem \ref{theo:deg_nf} we will assume that
\smallskip

(1) $\omega\in\mR^n$ satisfies the Bruno condition,

(2) $\widehat H\in\calF_\diamond\cap\calA$,

(3) the normal form of the Hamiltonian $H_2 + \widehat H$ is $H_2+N_\diamond$ $N_\diamond=O_r(\bz)$, $r\ge 3$.
\smallskip

Then by Lemma \ref{lem:sublinear-Bruno} there exists a sublinear convex negative sequence $\{b_s\}$ such that
\begin{equation}
\label{ab}
  a_j = \exp\big(b_{2^{j+1}+2} - 2 b_{2^j+2}\big)
\end{equation}
where $\{a_j\}$ is the following Bruno sequence
\begin{equation}
\label{a=Omega}
  a_j = n 2^{2n+j} (2^j+2)^{2n+1} \Omega_{2^j+2}, \qquad  j\in\mN.
\end{equation}

By Lemma \ref{lem:maj_sublinear} there exist $\widehat c_0$ and $\widehat\alpha_0$ such that
$$
  \widehat H = \sum_{|\bk|\ge 3} \widehat H_\bk \bz^\bk, \qquad
               |\widehat H_\bk|
           \le \widehat c_0 e^{b_{|\bk|} + \widehat\alpha_0 |\bk|}.
$$
Then, choosing $\alpha_0\ge\widehat\alpha_0$, we have:
\begin{equation}
\label{Hk<ce}
  |\widehat H_\bk| \le c_0 e^{b_{|\bk|}+\alpha_0 |\bk|}, \qquad
  |\bk| \ge 3,
\end{equation}
where $c_0 = \widehat c_0 e^{3(\widehat\alpha_0 - \alpha_0)}$. This means that choosing $\alpha_0$ sufficiently big, we may assume that
\begin{equation}
\label{calpha}
  c_0 e^{2\alpha_0} \le 1/16, \quad
  n e^{2\alpha_0} \ge 2.
\end{equation}

\begin{theo}
\label{theo:BR}
Suppose that $\omega$ and $\widehat H$ satisfy the above conditions (1)--(3) and (\ref{ab})--(\ref{calpha}).

Then an analytic canonical transformation of coordinates
\begin{equation}
\label{rho*}
  \thet : D_{\rho_*} \to D_{\rho_0}, \quad
  \bw\mapsto \bz = \thet(\bw), \qquad
  \rho_0 = e^{-\alpha_0}, \quad
  \rho_* = \rho_0 e^{-1/2}.
\end{equation}
reduces the Hamiltonian $H_2 + \widehat H$ to the form
$H_2\circ\thet + \widehat H\circ\thet = H_2 + G$, where
\begin{equation}
\label{beta>alpha}
   G(\bw)  =  \sum_{|\bk| \ge r} G_\bk \bw^\bk = O_r(\bw), \quad
  |G_\bk| \le c_0 e^{b_{|\bk|}+\beta |\bk|}, \qquad
   \beta  \le \alpha_0 + 1/2.
\end{equation}
\end{theo}

{\bf Remark}. Normalization we perform in the proof of Theorem \ref{theo:BR} deals with a finite number of ``small divisors'' $\langle\omega,\bk'\rangle$. Hence the assumption that $\omega$ satisfies the Bruno condition looks too strong.\footnote
{For $|\bk|>r$ we do not need any conditions on small divisors.}
We use this assumption to obtain uniform in $r$ estimates (\ref{rho*})--(\ref{beta>alpha}) for $\rho_*$ and $G_\bk$.

\begin{cor}
\label{cor:BR} (The Bruno-Russmann theorem).
Suppose that

(a) $\omega\in\mR^n$ satisfies the Bruno condition,

(b) $\widehat H \in\calA\cap\calF_\diamond$,

(c) the normal form of the Hamiltonian $H_2 + \widehat H$ equals $H_2$.
\smallskip

Then an analytic canonical transformation of coordinates $\bw\mapsto\bz=\thet(\bw)$ reduces the Hamiltonian $H_2+\widehat H$ to the normal form $H_2\circ\thet + \widehat H\circ\thet = H_2$.
\end{cor}

To prove Corollary \ref{cor:BR} it is sufficient to pass in Theorem \ref{theo:BR} to the limit
$r\to +\infty$.
\medskip

\subsection{Proof of Theorem \ref{theo:BR}}

The initial Hamiltonian satisfies estimate (\ref{Hk<ce}). We apply Lemma \ref{lem:low} inductively. Let
$$
  N = \max\{m\in\mZ : 2^{m-1} + 2 < r\}.
$$
Then $s$ in Lemma \ref{lem:low} takes the values
$$
  s_1,s_2,\ldots,s_N \qquad
  s_m = 2^{m-1} + 2.
$$
We obtain the sequences $\alpha_0,\alpha_1,\ldots,\alpha_N$, $\eps_1,\ldots,\eps_N$, $\Lambda_1,\ldots,\Lambda_N$, which estimate the Hamiltonian on $m$-th step:
\begin{equation}
\label{ael}
  \alpha_{m+1} = \alpha_m + \eps_{m+1}, \quad
   \eps_m   = c_0 \Lambda_m \Omega_{s_m}, \quad
  \Lambda_m = \frac{n e^{2\alpha_{m-1}}}{2} (2s_m)^{2n+1} e^{2b_{s_m} - b_{s_{m+1}}}.
\end{equation}
By (\ref{ab})--(\ref{a=Omega}) with $j=m-1$
\begin{equation}
\label{=1}
  e^{2b_{s_m} - b_{s_{m+1}}} n 2^{m-2} (2 s_m)^{2n+1} \Omega_{s_m} = 1.
\end{equation}
Equation (\ref{=1}) combined with (\ref{ael}) implies
\begin{equation}
\label{epsilon}
\!\!\!
    \eps_m
  = \frac{c_0 n e^{2\alpha_{m-1}}}{2} (2s_m)^{2n+1} e^{2b_{s_m} - b_{s_{m+1}}} \Omega_{s_m}
  = \frac{c_0 e^{2\alpha_{m-1}}}{2^{m-1}}
  = \frac{2c_0 e^{2\alpha_0 + 2(\eps_0+\ldots+\eps_{m-1})}}{2^m}.
\end{equation}

Now we prove by induction that for $c_0 e^{2\alpha_0}\le 1/16$ (see (\ref{calpha})) we have $\eps_m\le 2^{-m-2}$. Indeed, for $m=0$ the estimate holds: by (\ref{epsilon}) $\eps_0 = 2c_0^2 e^{2\alpha_0} \le 1/8$. Suppose it is valid for $m\le m_0$. Then
$$
      \eps_{m_0+1}
  \le \frac18 \frac{e^{2^{-2}+\ldots+2^{m-3}}}{2^{m_0+1}}
   <  \frac{e^{1/2}}{2^{m_0+2}}
   <  2^{-m_0-3}.
$$
Putting $\beta=\alpha_N$, we obtain:
$\beta - \alpha_0 \le \sum_{m=1}^\infty \eps_m \le 1/2$.
This implies (\ref{beta>alpha}).

To estimate $\rho_*$, we define the sequence $\thet_1,\ldots,\thet_N$, where
$\thet_m : D_{\rho_{m-1}}\to D_{\rho_m}$ is the map from Lemma \ref{lem:low} on the $m$-th step. We also define
\begin{equation}
\label{rho0rhoN}
  \rho_0,\ldots,\rho_N, \qquad
  \rho_0 = e^{-\alpha_0}, \quad
  \rho_N = \rho_* .
\end{equation}
This implies that the map $\thet=\thet_1\circ\ldots\circ\thet_N : D_{\rho_N}\to D_{\rho_0}$ is well-defined.

Inequality (\ref{rhorho'}) on $m$-th step reads
\begin{equation}
\label{rho<rho}
  \rho_{m+1} \le \rho_m - \frac{\eps_{m+1}}{n s e^\alpha_m} (\rho_m e^{\alpha_m})^{s_{m+1}-1}.
\end{equation}

Now we show that $\rho_m$ $(m\ge 1)$, determined by the equation $\rho_m e^{\alpha_m}=1$, satisfies (\ref{rho<rho}). Indeed, $\rho_0 e^{\alpha_0}=1$ by (\ref{rho0rhoN}). Assuming $\rho_m e^{\alpha_m}=1$ for some $m\ge0$, we have:
$$
    \bigg( \rho_m - \frac{\eps_{m+1}}{ns e^{\alpha_m}} (\rho_m e^{\alpha_m})^{s_{m+1}}
    \bigg)
  = e^{\eps_{m+1}} - \frac{\eps_{m+1}}{ns e^{\alpha_m}} e^{\alpha_{m+1}}
  = e^{\eps_{m+1}} \Big( 1 - \frac{\eps_{m+1}}{ns} \Big).
$$
The last quantity exceeds 1 because $\eps_m\le\eps_0\le 1/8$. Hence $\rho_{m+1}=e^{-\alpha_{m+1}}$ satisfies (\ref{rho<rho}).

Therefore $\rho_* =\rho_N = \rho_0 e^{-\sum_1^N \eps_m} \ge \rho_0 e^{-1/2}$.   \qed

\subsection{Proof of Theorem \ref{theo:deg_nf}}
\label{sec:proof_th}

Let $\{b_s\}$ be a sublinear convex negative sequence, satisfying (\ref{ab}) with the Bruno sequence (\ref{a=Omega}). Then by Theorem \ref{theo:BR} there exists an analytic canonical transformation of coordinates which reduces the Hamiltonian $H_2 + \widehat H$ to the form $H_2 + G$, where the function $G$ satisfies (\ref{beta>alpha}). Since $b_{|\bk|}$ are negative, we have:
$$
  |G_\bk| \le c_0 e^{\beta |\bk|}, \qquad
  |\bk| \ge r.
$$

Then we start another continuous averaging procedure with the operator $\xi=\xi_*$. We repeat the proof of Theorem \ref{theo:analytic} until equation (\ref{Bur_sol}). In our case $f'=O_{r-1}(\bz)$. By Lemma \ref{lem:maj1}
$$
  f'(\zeta) \ll \frac{a\zeta^{r-1}}{b-\zeta}, \qquad
  a = a(c_0,\beta), \quad
  b = b(c_0,\beta).
$$
Putting $\tau=8n\delta$, we obtain
\begin{equation}
\label{G=frac}
  G = \frac{a (\zeta + \tau G)^{r-1}}{b - (\zeta+\tau G)}.
\end{equation}

To estimate the function $G=G(\zeta)$, we have to use a quantitative version of the implicit function theorem (Lemma \ref{lem:Dieu} below). First, we introduce the function
$\widetilde G(\zeta) = \zeta + \tau G(\zeta)$. It satisfies the equation, obtained from (\ref{G=frac}):
\begin{equation}
\label{tildeG}
  \zeta = \widetilde G + \ph(\widetilde G), \qquad
  \ph(\widetilde G) = -\frac{a\tau\widetilde G^{r-1}}{b - \widetilde G}.
\end{equation}

\begin{lem}
\label{lem:Dieu} (see \cite{Dieu}).
Let $x=y+\ph(y)$ in the complex ball $\{y\in\mC : |y|\le 6\rho\}$, where the function $\ph$ is analytic and $|\ph|\le\rho/2$. Then there is a unique analytic function $\psi$ defined for
$|x|\le\rho$ such that
$$
  y = x + \psi(x), \quad
  |\psi| \le |\ph| \quad
  \mbox{for } |x| \le \rho.
$$
\end{lem}

Applying  Lemma \ref{lem:Dieu} to equation (\ref{tildeG}) with $x=\zeta$ and $y=\widetilde G$, we take
$$
  \rho = \min\Big\{ \frac{b}{12}, \frac16 \Big(\frac b{24 a\tau}\Big)^{1/(r-2)} \Big\}.
$$
Then in the ball $|\widetilde G| \le 6\rho$ we have
$$
      |\ph|
  \le \frac{a\tau (6\rho)^{r-1}}{b - 6\rho}
  \le \frac{2a\tau (6\rho)^{r-1}}{b}
  \le \frac \rho 2.
$$
This implies that
$\widetilde G = \zeta + \tau G = \zeta + \psi(\zeta)$,
where $|\psi| \le \rho/2$  for $|\zeta| \le \rho$.

If $\tau = 8n\delta$ is sufficiently large then
$$
  \rho = \frac16 \Big(\frac b{24 a\tau}\Big)^{1/(r-2)}
      \ge \frac{c_z}{1 + \delta^{1/(r-2)}}
$$
for some constant $c_z$. Hence in the domain $\{|\zeta|\le \rho\}$
$$
     |G|
  =  \Big|\frac{\widetilde G - \zeta}{\tau}\Big|
 \le \frac{\rho}{2\tau}
 \le \frac{C}{1 + \delta^{1 + 1/(r-2)}}
$$
for some constant $C$.
The function $\int_0^\zeta G(\widehat\zeta)\, d\widehat\zeta$ is a majorant for $|\calH|$.
This implies $|\calH| \le |\rho G| \le C_H / (1 +\delta^{1+2/(r-2)})$.
\qed

Generically $\omega$ does not admit resonances of order less than 4 (this includes the case of nonresonant $\omega$). Then $r = 4$ and diameter of the analyticity domain for $\calH$ is of order $\delta^{-1/2}$. This estimate is better than $\delta^{-1}$ declared in \cite{Tre23}.

\section{Technical part}
\label{sec:tech}

\subsection{Majorants}
\label{sec:maj}

\begin{lem}
\label{lem:maj_sublinear}
Suppose $F = \sum_{\bk\in\mZ_+^{2n}} F_\bk \bz^\bk \in\calA^\rho$, and the sequence $\{b_j\}_{j\in\mZ_+}$ is sublinear. Then for any $\alpha > - \ln\rho$ there exists $c>0$ such that
\begin{equation}
\label{|Fk|}
  |F_\bk| \le c e^{b_{|\bk|} + \alpha |\bk|}.
\end{equation}
\end{lem}

{\it Proof}. By Lemma \ref{lem:Hkk} we have $|F_\bk|\le c_F \rho^{-|\bk|}$, $c_F = \|F\|_\rho$. Hence in (\ref{|Fk|}) we have to choose $c$ such that
\begin{equation}
\label{ccF}
  c_F e^{-b_q - (\alpha + \ln\rho) q} < c\quad
  \mbox{for any }\;  q\in\mZ_+.
\end{equation}
For any $\alpha > - \ln\rho$ the function $q\mapsto e^{-b_q - (\alpha + \ln\rho) q}$, $q\in\mZ_+$ is bounded because $b_q$ is sublinear.  \qed

For any $F,\bF\in\calF$ we say that $F\ll\bF$ iff for their Taylor coefficients we have the inequalities $|F_\bk| \le \bF_\bk$, $\bk\in\mZ_+^{2n}$. In this case we say that $\bF$ is a majorant for $F$.

\begin{lem}
\label{lem:maj}
Suppose $F\ll\bF$ and $\hat F\ll\hat\bF$. Then

{\bf 1}. $F+\hat F\ll\bF+\hat\bF$, $F\hat F\ll \bF\hat\bF$,
$\partial_{z_s} F\ll\partial_{z_s}\bF$, and
         $\partial_{\overline z_s} F\ll\partial_{\overline z_s}\bF$,\, $s=1,\ldots,n$.

{\bf 2}. If $F$ and $\bF$ depend on the parameter $\delta\in [\delta_1,\delta_2]$ then
$$
      \int_{\delta_1}^{\delta_2} F\, d\delta
  \ll \int_{\delta_1}^{\delta_2} \bF\, d\delta .
$$
\end{lem}

We skip an obvious proof. \qed
\smallskip

\begin{lem}
\label{lem:maj1}
Suppose $F\in\calA^\rho$, $F = \sum_{|\bk|\ge s} F_\bk \bz^\bk = O_s(\bz)$, and
$|F_\bk| \le a\rho^{s - |\bk|}$. Then
\begin{equation}
\label{geom}
  F\ll \frac{a\rho\zeta^s}{\rho - \zeta}, \qquad
  \zeta = z_1 + \ldots + z_n + \overline z_1 + \ldots + \overline z_n.
\end{equation}
\end{lem}

{\it Proof}. Note that $\sum_{|\bk|=j} \bz^\bk \ll \zeta^j$ for any $j\in\mZ_+$. Then
$$
     F
 \ll \sum_{|\bk|\ge s} a\rho^{s - |\bk|} \bz^\bk
 \ll a\rho^s \sum_{j=s}^\infty \frac{\zeta^j}{\rho^j}
  =  \frac{a\rho\zeta^s}{\rho - \zeta}.
$$
\qed

\begin{lem}
\label{lem:f'}
For any $\rho>0$
$$
      \partial_\zeta \frac{\rho\zeta^3}{\rho - \zeta}
  \ll \frac{2\rho\zeta^2}{\rho/2 - \zeta}.
$$
\end{lem}

{\it Proof}. Since
$$
    \partial_\zeta \frac{\rho\zeta^3}{\rho - \zeta}
  = \frac{3\rho\zeta^2}{\rho - \zeta} + \frac{\rho\zeta^3}{(\rho - \zeta)^2},
$$
the lemma follows from two estimates:
\begin{eqnarray*}
       \frac{3\rho\zeta^2}{\rho - \zeta}
 &\ll& \frac{3\rho\zeta^2}{2(\rho/2 - \zeta)}, \\
       \frac{\rho\zeta}{(\rho - \zeta)^2}
 &=&   \sum_{k=0}^\infty k \Big(\frac\zeta\rho\Big)^k
  =    \sum_{k=0}^\infty \frac{k}{2^k} \Big(\frac\zeta{\rho/2}\Big)^k
 \ll   \sum_{k=0}^\infty \frac12 \Big(\frac\zeta{\rho/2}\Big)^k
  =    \frac{\rho}{2(\rho/2 - \zeta)^2}.
\end{eqnarray*}
\qed

\subsection{Majorant principle}
\label{sec:MP}

We use majorant method to obtain estimates for solutions of initial value problems (IVP) in $\calF$.

As an example consider the IVP
\begin{equation}
\label{F'=Phi}
  \partial_\delta F = \Phi(F,\delta), \qquad
  F|_{\delta=0} = \widehat F.
\end{equation}
Here $F\in\calF$ depends on the parameter $\delta$ and $\Phi$ is a map from $\calF\times\mR_+$ to $\calF$.

\begin{dfn}
\label{dfn:pr}
IVP (\ref{F'=Phi}) is said to be power regular if for any $\widehat F\in\calF$ equation (\ref{F'=Phi}) has a unique solution $F = F(\bz,\delta)\in\calF$ for all $\delta>0$
\end{dfn}

We associate with (\ref{F'=Phi}) the so-called majorant system
\begin{equation}
\label{bF'=Psi}
  \partial_\delta \bF = \Psi(\bF,\delta), \qquad
  \bF|_{\delta=0} = \widehat\bF.
\end{equation}

We put $\Phi_\bk = p_\bk\circ\Phi$ and $\Psi_\bk = p_\bk\circ\Psi$.

\begin{dfn}
\label{dfn:maj}
IVP (\ref{bF'=Psi}) is said to be a majorant IVP for (\ref{F'=Phi}) if the following two properties hold:

(a) $\widehat F\ll \widehat\bF$.

(b) For any $F\ll\bF$, $\bk\in\mZ_+^{2n}$, and $\delta \ge 0$ we have: $\Phi_\bk(F,\delta) \ll \Psi_\bk(\bF,\delta)$.
\end{dfn}

{\bf Majorant principle}. {\it Suppose the IVP (\ref{F'=Phi}) is power regular. Suppose also that there exists a solution $\bF = \bF(\cdot,\delta)\in\calA$ of (\ref{bF'=Psi}) on the interval
$\delta\in [0,\delta_0]$. Then (\ref{F'=Phi}) has a unique analytic solution $F$ on $[0,\delta_0]$ and $F(\cdot,\delta) \ll \bF(\cdot,\delta)$.
}
\medskip

{\bf Remarks}. 1. Definitions \ref{dfn:pr} and \ref{dfn:maj} as well as the Majorant principle obviously extend to systems of equations, where $F,\widehat F\in\calF^m$ and $\Phi:\calF^m\times\mR_+\to\calF$.

2. One may replace the first equation (\ref{bF'=Psi}) by the inequality
$\partial_\delta \bF \gg \Psi(\bF,\delta)$.

\medskip

The majorant principle presented here differs from the majorant argument used since Cauchy times. Traditionally the evolution variable (in our case $\delta$) is regarded complex as well and Taylor expansions in it are used. In our approach this variable is a real parameter in both exact solution and a majorant. Due to this we are able to obtain majorant estimates for solutions of (\ref{F'=Phi}) on large (even infinite) intervals of $\delta$.
\smallskip

\begin{theo}
\label{theo:maj_nil}
Suppose both systems (\ref{F'=Phi}) and (\ref{bF'=Psi}) have nilpotent structure (see Definition \ref{dfn:nil}). Then Majorant principle holds true.
\end{theo}

We expect that Majorant principle is valid in a much wider generality. But in this paper we are only interested in the case of systems having nilpotent form.
\smallskip

{\it Proof of Theorem \ref{theo:maj_nil}}.
Let $\bk^0$ be an index with minimal possible degree $|\bk^0|$. For example, in the system (\ref{aver4}) $|\bk^0|=3$. Nilpotent form of (\ref{F'=Phi}) implies that
$$
  0 = \partial_\delta F_{\bk^0} \ll \partial_\delta \bF_{\bk^0}.
$$
Hence $F_{\bk^0}(\delta)\ll\bF_{\bk^0}(\delta)$ for $\delta\ge 0$.

We proceed by induction in $|\bk|$. Suppose $F_{\bk}(\delta)\ll\bF_{\bk}(\delta)$, $\delta\ge 0$ provided $|\bk|< K$. For any $\bk$ such that $|\bk|=K$ we have by induction assumption and item (b) of Definition \ref{dfn:maj}:
$$
    \partial_\delta (\bF_\bk - F_\bk)
  = \Psi_\bk(\bF(\cdot,\delta),\delta) - \Phi_\bk(\bF(\cdot,\delta),\delta)
 \gg 0.
$$
Therefore
$$
    \bF_\bk(\delta)
  = \widehat\bF_\bk(\delta)
   + \int_0^\delta \Psi_\bk(\bF(\cdot,\lambda),\lambda)\, d\lambda
 \gg F_\bk(\delta)
  = \widehat F_\bk(\delta)
   + \int_0^\delta \Phi_\bk(F(\cdot,\lambda),\lambda)\, d\lambda.
$$
Here we used that arguments of $\Psi_\bk$ and $\Phi_\bk$ are known by the induction assumption.

This majorant inequality makes sense if the left-hand side is defined i.e., for any
$\delta\in [0,\delta_0]$.  \qed

\subsection{On convex sequences}

\begin{lem}
\label{lem:conv_seq}
Suppose the sequence $\{b_j\}_{j\in\mN}$ is convex. Then

(1) for any $1\le m < k < l$
\begin{equation}
\label{b+b+b}
   (l-k)b_m + (m-l)b_k + (k-m)b_l \ge 0.
\end{equation}

(2) for any $1\le m < k \le l$
\begin{equation}
\label{b+b}
   b_k + b_l \le b_{k-m} + b_{l+m},
\end{equation}
\end{lem}

{\it Proof}. (1) First, consider the case $l-k=1$. If $k-m=1$ then (\ref{b+b+b}) coincides with definition of convexity. We use induction in $k-m$. Suppose (\ref{b+b+b}) holds for $k-m=q$. Then inequality (\ref{b+b+b}) is obtained if we add the inequalities
\begin{eqnarray*}
&&  (l-k) b_m + (m-l+1) b_{k-1} + (k-1-m) b_{l-1} \ge 0, \\
&&\qquad\qquad\;     -(m-l+1) b_{k-1} + 2(m-l+1) b_{l-1} - (m-l+1) b_l \ge 0
\end{eqnarray*}
which hold by the induction assumption.

The case $l-k>1$ follows from the case $l-k=1$ also by induction.

(2) Equation (\ref{b+b+b}) implies the inequalities
\begin{eqnarray*}
       (l-k) b_{k-m} + (k-m-l) b_k + m b_l
 &\ge& 0, \\
       m b_k + (k-l-m) b_l + (l-k) b_{l+m}
 &\ge& 0
\end{eqnarray*}
from which (\ref{b+b}) follows. \qed

\begin{lem}
\label{lem:sublinear-Bruno}
{\bf 1}. For any convex sublinear sequence $\{b_s\}$ the sequence $\{a_j\}$, satisfying (\ref{ab})
is Bruno.

{\bf 2}. For any Bruno sequence $\{a_j\}$ there exists a sublinear convex negative non-increasing sequence $\{b_s\}$, satisfying (\ref{ab}).
\end{lem}

{\it Proof}. {\bf 1}. Suppose $\{b_s\}$ is convex and sublinear, and $\{a_j\}$ satisfies (\ref{ab}). First, note that
$$
  \ln a_j - \ln a_{j-1} = - b_{2^{j+1}+2} + 3b_{2^j+2} - 2b_{2^{j-1}+2}.
$$
The latter expression is nonnegative by convexity of $\{b_s\}$ (see Lemma \ref{lem:conv_seq}). Hence, $a_j$ is non-decreasing.

Equation (\ref{ab}) implies
$$
    \frac12 \sum_{j=0}^{J-1} 2^{-j} \ln a_j
  = \sum_{j=0}^{J-1} \Big( 2^{-j-1} b_{2^{j+1}+2} - 2^{-j} b_{2^j+2} \Big)
  = 2^{-J} b_{2^J+2} - b_3.
$$
Since $\{b_j\}$ is sublinear, $2^{-J} b_{2^J+2} \to 0$ as $J\to\infty$. Therefore $\{a_j\}$ is a Bruno sequence.

{\bf 2}. Suppose $\{a_j\}$ is a Bruno sequence. We compute $\{b_{2^j+2}\}$ from the equation
$$
  \frac{b_{2^J+2}}{2^J} = b_3 + \frac12 \sum_{j=1}^{J-1} \frac{\ln a_j}{2^j}
$$
which implies (\ref{ab}). To have sublinearity, we need $\lim_{J\to +\infty} 2^{-J} b_{2^J+2} = 0$. Hence, we take $b_3 = - \frac12 \sum_{j=1}^\infty 2^{-j} \ln a_j$. Such a choice of $b_3$ implies that $b_{2^J+2} < 0$, $J=0,1,\ldots$

The equation
$$
    b_{2^J+2}
  = - 2^{J-1} \sum_{j=J}^\infty  \frac{\ln a_j}{2^j}
  = - \frac12 \sum_{j=0}^\infty  \frac{\ln a_{J+j}}{2^j}
$$
implies
$$
    b_{2^{J+1}+2} - b_{2^J+2}
  = \frac12 \ln a_J - \frac12 \sum_{j=1}^\infty  \frac{\ln a_{J+j}}{2^j}
 \le \frac12 \ln a_J \Big( 1 - \sum_{j=1}^\infty 2^{-j} \Big)
  =  0.
$$
Hence the sequence $\{b^{2^j+2}\}$ does not increase.

For any integer $s\in (2^j+2,2^{j+1}+2)$ we define $b_s$ by linear interpolation:
$$
    b_s
  = \frac{s - 2^j - 2}{2^j} b_{2^j+2} + \frac{2^{j+1} + 2 - s}{2^j} b_{2^{j+1}+2}.
$$
Such $\{b_s\}$ is obviously sublinear.

We check convexity of $\{b_s\}$ first, on the subsequence $\{b_{2^j+2}\}_{j\in\mZ_+}$. We have to prove the inequality
\begin{equation}
\label{ineq}
  2^J b_{2^{J-1}+2} - (2^J + 2^{J-1}) b_{2^J+2} + 2^{J-1} b_{2^{J+1}+2} \ge 0.
\end{equation}
By (\ref{ab}) the left-hand side of (\ref{ineq}) equals $2^{J-1}(\ln a_J - \ln a_{J-1})$. For any
$J\ge 1$ this expression is nonnegative because $\{a_j\}$ is non-decreasing.

If $b_s$, $s\ne 2^j+2$ are determined by linear interpolation, the convexity, negativity, and monotonicity remain true for all the sequence $\{b_s\}_{s\in\mN}$.  \qed

\end{document}